\documentstyle[11pt,amssymb,amstex]{amsart}

\numberwithin{equation}{section}

\newtheorem{theorem}{Theorem}[section]
\newtheorem{proposition}[theorem]{Proposition}
\newtheorem{lemma}[theorem]{Lemma}

\newtheorem{definition}[theorem]{Definition}
\newtheorem{remark}[theorem]{Remark}

\def\bn{{\mathbb N}}

\def\a{\alpha}

\def\<{\langle}
\def\>{\rangle}

\def\1{\mathbf{1}}

\def\xb{\mathbf{x}}
\def\yb{\mathbf{y}}

\begin{document}

\title[Quasi $I-$nonexpansive mappins]
{Weak and strong convergence of an implicit iterative process with
errors for a finite family of asymptotically quasi $I-$nonexpansive
mappings in Banach space}

\author{Farrukh Mukhamedov}
\address{Farrukh Mukhamedov\\
Department of Computational \& Theoretical Sciences \\
Faculty of Sciences, International Islamic University Malaysia\\
P.O. Box, 141, 25710, Kuantan\\
Pahang, Malaysia} \email{{\tt far75m@@yandex.ru}}

\author{Mansoor Saburov}
\address{Mansoor Saburov\\
Department of Computational \& Theoretical Sciences \\
Faculty of Science, International Islamic University Malaysia\\
P.O. Box, 141, 25710, Kuantan\\
Pahang, Malaysia} \email{{\tt msaburov@@gmail.com}}

\begin{abstract}
In this paper we prove the weak and strong convergence of the
implicit iterative process with errors to a common fixed point of a finite family $\{T_j\}_{i=1}^N$ of  asymptotically quasi $I_j-$nonexpansive mappings as well as a family of $\{I_j\}_{j=1}^N$ of asymptotically quasi nonexpansive mappings in the framework
of Banach spaces.
\vskip 0.3cm \noindent
{\it
Mathematics Subject Classification}: 46B20; 47H09; 47H10\\
{\it Key words}: Implicit iteration process with errors; a finite
family of asymptotically quasi $I-$nonexpansive mapping; asymptotically
quasi-nonexpansive mapping; a common fixed point; Banach space.
\end{abstract}

\maketitle

\section{Introduction}
Let $K$ be a nonempty subset of a real normed linear space $X$ and $T:K\to K$ be a mapping.
 Denote by $F(T)$ the set of fixed points of $T$, that is, $F(T) =\{x\in K: Tx = x\}$. Throughout this paper,
 we always assume that $F(T)\neq\emptyset$. Now let us recall some known definitions
\begin{definition}
A mapping $T:K\to K$ is said to be:
\begin{itemize}
  \item[(i)] nonexpansive, if $\|Tx-Ty\|\le\|x-y\|$
  for all $x,y\in K$;
  \item[(ii)]  asymptotically nonexpansive,
  if there exists a sequence $\{\lambda_{n}\}\subset[1,\infty)$ with
  $\lim\limits_{n\to\infty}\lambda_{n}=1$ such that
  $\|T^nx-T^ny\|\le\lambda_n\|x-y\|$ for all $x,y\in K$ and $n\in\bn$;
  \item[(iii)] quasi-nonexpansive, if
  $\|Tx-p\|\le\|x-p\|$ for all $x\in K,\ p\in F(T)$;
  \item[(iv)] asymptotically quasi-nonexpansive,
  if there exists a sequence $\{\mu_n\}\subset [1,\infty)$
  with $\lim\limits_{n\to\infty}\mu_n=1$ such that
  $\|T^nx-p\|\le\mu_n\|x-p\|$ for all $x\in K,\ p\in F(T)$ and $n\in\bn$.
  \end{itemize}
\end{definition}

Note that from the above definitions, it follows that a nonexpansive mapping must be asymptotically
nonexpansive, and an asymptotically nonexpansive mapping must be asymptotically
quasi-nonexpansive, but the converse does not hold (see \cite{GK2}).


If $K$ is a closed nonempty subset of a Banach space and $T:K\to K$
is nonexpansive, then it is known that $T$ may not have a fixed point
(unlike the case if $T$ is a strict contraction), and even when it
has, the sequence $\{x_n\}$ defined by $x_{n+1} = Tx_n $ (the
so-called \emph{Picard sequence}) may fail to converge to such a
fixed point.

In \cite{[Browder65]}-\cite{[Browder67]} Browder  studied the iterative
construction for fixed points of nonexpansive mappings on closed and
convex subsets of a Hilbert space. Note that for the past 30 years or so, the study of the iterative processes
for the approximation of fixed points of nonexpansive mappings and
fixed points of some of their generalizations have been flourishing
areas of research for many mathematicians (see for more details \cite{GK2},\cite{Chidumi}).

In \cite{[DaizMetcalf67]} Diaz and Metcalf
studied quasi-nonexpansive mappings in Banach spaces. Ghosh and
Debnath \cite{[GhoshDebnath]} established a necessary and sufficient
condition for convergence of the Ishikawa iterates of a
quasi-nonexpansive mapping on a closed convex subset of a Banach
space. The iterative approximation problems for nonexpansive mapping, asymptotically
nonexpansivemapping and asymptotically quasi-nonexpansive mapping were
studied extensively by  Goebel and Kirk \cite{[GeobelKrik]}, Liu \cite{liu}, Wittmann \cite{[Wittmann]}, Reich
\cite{R}, Gornicki \cite{[Gornicki]}, Shu \cite{[Shu]} Shoji and Takahashi \cite{ST}, Tan and Xu \cite{[TanXu]} et al. in the settings of Hilbert spaces and
uniformly convex Banach spaces.

There are many methods for approximating fixed points of a
nonexpansive mapping. Xu and Ori \cite{XO} introduced implicit
iteration process to approximate a common fixed point of a finite
family of nonexpansive mappings in a Hilbert space. Namely, let $H$
be a Hilbet space, $K$ be a nonempty closed convex subset of $H$ and
$\{T_j\}_{j=1}^N: K\to K$ be nonexpansive mappings. Then Xu and
Ori's implicit iteration process $\{x_n\}$ is defined by
\begin{eqnarray}\label{xo}
\left\{ \begin{array}{cc}
x_0\in K,\\
x_n = (1-\a_n)x_{n-1}+ \a_nT_nx_n;\\
\end{array}
\right. \ \ n\geq 1
\end{eqnarray}
where $T_n = T_{n (mod N)}$, $\{\a_n\}$ is a real sequence in $(0,
1)$. They proved the weak convergence of the sequence $\{x_n\}$
defined by \eqref{xo} to a common fixed point $p\in F =\cap_{j=1}^N
F(T_j)$.

In 2003, Sun \cite{[Sun]} introduced the following implicit iterative sequence
$\{x_n\}$
\begin{eqnarray}\label{sun}
\left\{ \begin{array}{cc}
x_0\in K,\\
x_n = (1-\a_n)x_{n-1}+ \a_nT_{j(n)}^{k(n)}x_n;\\
\end{array}
\right. \ \ n\geq 1
\end{eqnarray}
for a finite set of asymptotically quasi-nonexpansive self-mappings on a
bounded closed convex subset $K$ of a Banach space $X$ with $\{\a_n\}$ a sequence
in $(0, 1)$.
Where $n =(k(n)-1)N + j(n)$, $j(n)\in\{1, 2,\dots,N\}$, and proved the strong convergence of the
sequence $\{x_n\}$ defined by \eqref{sun} to a common fixed point $p\in F =\cap_{j=1}^N F(T_j)$.

There many papers devoted to the implicit iteration process for a finite  family
of asymptotically nonexpansive mappings, asymptotically quasi-expansive mappings in Banach spaces (see
for example \cite{CTL,CS,G,GL,liu,ZC}).

On the other hand, there are many concepts which generalize a notion of nonexpansive mapping. One of such concepts is
$I$-nonexpansivity of a mapping $T$(\cite{Shah}). Let us recall some notions.

\begin{definition}
Let $T:K\to K$, $I:K\to K$ be two mappings of a nonempty subset $K$ of a real normed linear space $X$.
Then $T$ is said to be:
\begin{itemize}
\item[(i)] {\it $I-$nonexpansive}, if $\|Tx-Ty\|\le\|Ix-Iy\|$ for all $x,y\in K$;
\item[(i)] {\it asymptotically $I-$nonexpansive}, if there exists a sequence $\{\lambda_{n}\}\subset[1,\infty)$ with   $\lim\limits_{n\to\infty}\lambda_{n}=1$ such that
  $\|T^nx-T^ny\|\le\lambda_n\|I^nx-I^ny\|$ for all $x,y\in K$ and $n\ge 1$;
\item[(iii)] {\it asymptotically quasi $I-$nonexpansive}
  mapping, if there exists a sequence $\{\mu_n\}\subset [1,\infty)$
  with $\lim\limits_{n\to\infty}\mu_n=1$ such that
  $\|T^nx-p\|\le\mu_n\|I^nx-p\|$ for all $x\in K,\ p\in F(T)\cap F(I)$ and $n\ge 1.$
\end{itemize}
\end{definition}

\begin{remark} If $F(T)\cap F(I)\neq \emptyset$ then an asymptotically
$I-$nonexpansive mapping is asymptotically quasi $I-$nonexpansive.
But, there exists a nonlinear continuous asymptotically quasi
$I-$nonexpansive mappings which is asymptotically $I-$nonexpansive.
\end{remark}

Indeed, let us  consider the following example. Let
 $X=\ell_2$ and $K=\{{\xb}\in\ell_2:\ \|{\xb}\|\leq 1\}$. Define the following mappings:
\begin{eqnarray}\label{T}
&&
T(x_1,x_2,\dots,x_n,\dots)=(0,x_1^4,x_2^4,\dots,x_n^4,\dots),\\[2mm]
\label{I}
&&I(x_1,x_2,\dots,x_n,\dots)=(0,x_1^2,x_2^2,\dots,x_n^2,\dots).
\end{eqnarray}
One see that $F(T)=F(I)=(0,0,\dots,0,\dots)$. Therefore, from
$\sum\limits_{k=1}^\infty x_k^8\leq\sum\limits_{k=1}^\infty x_k^4$
whenever ${\xb}\in K$, using \eqref{T},\eqref{I} we obtain
$\|T{\xb}\|\leq\|I{\xb}\|$ for every ${\xb}\in K$. So, $T$ is quasi
$I$-expansive. But for ${\xb}_0=(1,0,\dots,0,\dots)$ and
${\yb}_0=(1/2,0,\dots,0,\dots)$ we have
$$
\|T({\xb}_0)-T({\yb}_0)\|=\frac{15}{16}, \ \ \
\|I({\xb}_0)-I({\yb}_0)\|=\frac{3}{4}
$$
which means that $T$ is not $I$-nonexpansive.

Note that in \cite{TG} the weakly convergence theorem for $I$-asymptotically quasi-nonexpansive
mapping defined in Hilbert space was proved.  Best approximation properties of $I$-nonexpansive
mappings were investigated in \cite{Shah}. In
\cite{KKJ} the weak convergence of three-step Noor iterative scheme for
an $I$-nonexpansive mappping in a Banach space has been established.

Very recently, in \cite{T} the weak and strong convergence of implicit iteration process to a common fixed point of a finite family of $I$-asymptotically nonexpansive mappings were studied. Let us describe the iteration scheme considered in \cite{T}. Let $K$ be a nonempty convex subset of a real Banach space $X$ and $\{T_j\}_{j=1}^N:K\to K$ be a finite family of asymptotically $I_j-$nonexpansive mappings, $\{I_j\}_{j=1}^N:K\to K$ be a finite
family of asymptotically nonexpansive mappings. Then the iteration process $\{x_n\}$ has been defined by
\begin{eqnarray}\label{implicitmap}
\left\{ \begin{array}{ccc}
x_1\in K,\\
          x_{n+1} = & (1-\alpha_n) x_{n}+\alpha_n I_{j(n)}^{k(n)}y_n,\\
          y_n = & (1-\beta_n) x_n+\beta_n
          T_{j(n)}^{k(n)}x_n,
        \end{array}\right. \ \ n\geq 1
\end{eqnarray}
here as before  $n=(k(n)-1)N+j(n)$, $j(n)\in\{1,2,\dots,N\}$, and  $\{\alpha_n\},$ $\{\beta_n\}$ are two sequences in $[0,1]$. From this formula one can easily see that the employed method, indeed, is not implicit iterative processes. The used process is some kind of modified Ishikawa iteration.

Therefore, in this paper we shall extend of the implicit iterative process with errors, defined in \cite{[Sun],fk}, to a family of $I$-asymptotically quasi-nonexpansive mappings defined on a uniformly convex Banach space. Namely,
let $K$ be a nonempty convex subset of a real Banach space $X$ and
$\{T_j\}_{j=1}^N:K\to K$ be a finite family of asymptotically quasi $I_j-$nonexpansive mappings, and
$\{I_j\}_{j=1}^N:K\to K$  be a family of asymptotically quasi-nonexpansive mappings. We consider the following implicit iterative scheme $\{x_n\}$ with errors:
\begin{eqnarray}\label{implicitmap}
\left\{ \begin{array}{ccc}
x_0\in K,\\
          x_n = & \alpha_n x_{n-1}+\beta_n T_{j(n)}^{k(n)}y_n+\gamma_nu_n\\
          y_n = & \widehat{\alpha}_n x_n+\widehat{\beta}_n
          I_{j(n)}^{k(n)}x_n+\widehat{\gamma}_nv_n
        \end{array}\right. \ \ n\geq 1
\end{eqnarray}
where $\{\alpha_n\},$ $\{\beta_n\},$ $\{\gamma_n\},$
$\{\widehat{\alpha}_n\},$ $\{\widehat{\beta}_n\}$,
$\{\widehat{\gamma}_n\}$ are six sequences in $[0,1]$ satisfying
$\alpha_n+\beta_n+\gamma_n=\widehat{\alpha}_n+\widehat{\beta}_n+\widehat{\gamma}_n=1$
for all $n\ge1$,  as well as
$\{u_n\}$, $\{v_n\}$ are bounded sequences in $K.$

In this paper we shall prove the weak and strong convergence of the
implicit iterative process \eqref{implicitmap} to a common fixed
points of $\{T_j\}_{j=1}^N$ and $\{I_j\}_{j=1}^N$. All results
presented here generalize and extend the corresponding main results
of \cite{[Sun]}, \cite{XO},\cite{fk},\cite{gou}.

\section{Preliminaries}

Throughout this paper, we always assume that $X$ is a real Banach
space. We denote $F(T)$ and $D(T)$ the set
of  fixed points and the domain of a mapping $T,$ respectively.
Recall that a Banach space $X$ is said to satisfy \emph{Opial
condition} \cite{[Opial]}, if for each sequence $\{x_{n}\}$ in $X,$
the condition that the sequence $x_n\to x$ weakly implies that
\begin{eqnarray}\label{Opialcondition}
\liminf\limits_{n\to\infty}\|x_n-x\|<
\liminf\limits_{n\to\infty}\|x_n-y\|
\end{eqnarray}
for all $y\in X$ with $y\neq x.$ It is well known that (see
\cite{[LamiDozo]}) inequality \eqref{Opialcondition} is equivalent
to
\begin{eqnarray*}
\limsup\limits_{n\to\infty}\|x_n-x\|<
\limsup\limits_{n\to\infty}\|x_n-y\|
\end{eqnarray*}

\begin{definition}
Let K be a closed subset of a real Banach space $X$ and $T:K\to K$
be a mapping.
\begin{itemize}
  \item[(i)]A mapping $T$ is said to be semi-closed (demi-closed) at zero, if for
  each bounded sequence $\{x_n\}$ in $K,$ the conditions $x_n$ converges weakly to $x\in K$ and $Tx_n$
  converges strongly to $0$ imply $Tx=0.$
  \item[(ii)] A mapping $T$ is said to be semi-compact, if for any
  bounded sequence $\{x_n\}$ in $K$ such that $\|x_n-Tx_n\|\to 0,
  (n\to\infty),$ then there exists a subsequence
  $\{x_{n_k}\}\subset\{x_n\}$ such that $x_{n_k}\to x^*\in K$
  strongly.
  \item[(iii)] $T$ is called a uniformly $L-$Lipschitzian mapping, if
  there exists a constant $L>0$ such that $\|T^nx-T^ny\|\le L\|x-y\|$ for all $x,y\in
  K$ and $n\ge 1.$
 \end{itemize}
\end{definition}

\begin{proposition}\label{commonLandlambda}
Let $K$ be a nonempty subset of a real Banach space $X,$
$\{T_j\}_{j=1}^N:K\to K$ and $\{I_j\}_{j=1}^N:K\to K$ be finite
familis of mappings.
\begin{itemize}
  \item[(i)] If $\{T_j\}_{j=1}^N$ is a finite family of
  asymptotically $I_j-$nonexpansive (resp. asymptotically quasi
  $I_j-$nonexpansive) mappings with sequences
  $\{\lambda^{(j)}_n\}\subset[1,\infty),$ then there exists a
  sequence $\{\lambda_n\}\subset[1,\infty)$ such that $\{T_j\}_{i=1}^N$ is a finite family of
  asymptotically $I_j-$nonexpansive (resp. asymptotically quasi
  $I_j-$nonexpansive) mappings with a common sequence
  $\{\lambda_n\}\subset[1,\infty).$
  \item[(ii)] If $\{T_j\}_{j=1}^N$ is a finite family of
  uniformly $L_j-$Lipschitzian mappings, then there exists a constant $L>0$ such that
  $\{T_j\}_{j=1}^N$ is a finite family of uniformly $L-$Lipschitzian
  mappings.
\end{itemize}
\end{proposition}
\begin{pf}
We shall prove (i) part of Proposition \ref{commonLandlambda}.
Analogously one can prove (ii) part.

Without any loss of generality, we may assume that $\{T_j\}_{j=1}^N$
is a finite family of asymptotically $I_j-$nonexpansive mappings
with sequences $\{\lambda^{(j)}_n\}\subset[1,\infty).$ Then we have
$$\|T^n_jx-T^n_jy\|\le\lambda^{(j)}_n\|I^n_jx-I^n_jy\|, \qquad \forall x,y\in K,\  \forall n\ge1,\  j=\overline{1,N}.$$
Let $\lambda_n=\max\limits_{j=\overline{1,N}}\lambda^{(j)}_n.$ Then
$\{\lambda_n\}\subset[1,\infty)$ with $\lambda_n\to 1,$ $n\to\infty$
and
$$\|T^n_jx-T^n_jy\|\le\lambda_n\|I^n_jx-I^n_jy\|, \ \  \forall n\ge1,$$
for all $x,y\in K,$ and each $j=\overline{1,N}.$ This means that
$\{T_j\}_{i=1}^N$ is a finite family of asymptotically
$I_j-$nonexpansive mappings with a common sequence
$\{\lambda_n\}\subset[1,\infty).$
\end{pf}

The following lemmas play an important role in proving our main
results.

\begin{lemma}[ see \cite{[Shu]}]\label{convexxnyn}
Let $X$ be a uniformly convex Banach space and $b,c$ be two
constants with $0<b<c<1.$ Suppose that $\{t_n\}$ is a sequence in
$[b,c]$ and $\{x_n\},$ $\{y_n\}$ are two sequences in $X$ such that
\begin{eqnarray*}
\lim\limits_{n\to\infty}\|t_nx_n+(1-t_n)y_n\|=d, \ \ \
\limsup\limits_{n\to\infty}\|x_n\|\le d, \ \ \
\limsup\limits_{n\to\infty}\|y_n\|\le d,
\end{eqnarray*}
holds some $d\ge 0.$ Then $\lim\limits_{n\to\infty}\|x_n-y_n\|=0.$
\end{lemma}

\begin{lemma}[see \cite{[TanXu]}]\label{convergean}
Let $\{a_n\},$ $\{b_n\}$ and $\{c_n\}$ be three sequences of
nonnegative real numbers with $\sum\limits_{n=1}^{\infty}b_n<\infty$
and $\sum\limits_{n=1}^{\infty}c_n<\infty.$ If the following
conditions is satisfied
$$a_{n+1}\le(1+b_n)a_n+c_n, \ \ \ n\ge 1,$$
then the limit $\lim\limits_{n\to\infty}a_n$ exists.
\end{lemma}

\section{Main results}

In this section we shall prove our main results concerning weak and
strong convergence of the sequence defined by \eqref{implicitmap}.
To formulate ones, we need some auxiliary results.

\begin{lemma}\label{limexistsxnminusp}
Let $X$ be a real Banach space and $K$ be a nonempty closed convex
subset of $X.$ Let $\{T_j\}_{j=1}^N:K\to K$ be a finite family of
asymptotically quasi $I_j-$nonexpansive mappings with a common
sequence $\{\lambda_n\}\subset[1,\infty)$ and $\{I_j\}_{j=1}^N:K\to
K$ be a finite family of asymptotically quasi-nonexpansive mappings
with a common sequence $\{\mu_n\}\subset[1,\infty)$ such that
$F=\bigcap\limits_{j=1}^N \left(F(T_j)\cap
F(I_j)\right)\neq\emptyset.$ Suppose $B^{*}=\sup\limits_{n}\beta_n,$
$\Lambda=\sup\limits_{n}\lambda_n\ge1,$ $M=\sup\limits_{n}\mu_n\ge1$
and $\{\alpha_n\},$ $\{\beta_n\},$ $\{\gamma_n\},$
$\{\widehat{\alpha}_n\},$ $\{\widehat{\beta}_n\},$
$\{\widehat{\gamma}_n\}$ are six sequences in $[0,1]$  which satisfy
the following conditions:
\begin{itemize}
  \item[(i)] $\alpha_n+\beta_n+\gamma_n=\widehat{\alpha}_n+\widehat{\beta}_n+\widehat{\gamma}_n=1, \ \ \ \forall n\ge1,$
  \item[(ii)] $\sum\limits_{n=1}^{\infty}(\lambda_n\mu_n-1)\beta_n<\infty,$
  \item[(iii)] $B^{*}<\dfrac{1}{\Lambda^2M^2},$
  \item[(iv)] $\sum\limits_{n=1}^\infty\gamma_n<\infty, \ \
  \sum\limits_{n=1}^\infty\widehat{\gamma}_n<\infty.$
\end{itemize}
Then for the implicit iterative sequence $\{x_n\}$ with errors
defined by \eqref{implicitmap} and for each $p\in F$ the limit
$\lim\limits_{n\to\infty}\|x_n-p\|$ exists.
\end{lemma}

\begin{pf}
Since $F=\bigcap\limits_{j=1}^N \left(F(T_j)\cap
F(I_j)\right)\neq\emptyset,$ for any given $p\in F,$ it follows from
\eqref{implicitmap} that
\begin{eqnarray*}
\|x_n-p\| &=& \|\alpha_n(x_{n-1}-p)+\beta_n(T^{k(n)}_{j(n)}y_n-p)+\gamma_n(u_n-p)\|\\
&\le& (1-\beta_n-\gamma_n)\|x_{n-1}-p\|+\beta_n\|T^{k(n)}_{j(n)}y_n-p\|+\gamma_n\|u_n-p\|\\
&\le& (1-\beta_n-\gamma_n)\|x_{n-1}-p\|+\beta_n\lambda_{k(n)}\|I^{k(n)}_{j(n)}y_n-p\|+\gamma_n\|u_n-p\|\\
&\le& (1-\beta_n-\gamma_n)\|x_{n-1}-p\|+\beta_n\lambda_{k(n)}\mu_{k(n)}\|y_n-p\|+\gamma_n\|u_n-p\|\\
&\le&
(1-\beta_n)\|x_{n-1}-p\|+\beta_n\lambda_{k(n)}\mu_{k(n)}\|y_n-p\|+\gamma_n\|u_n-p\|.
\end{eqnarray*}
Again using \eqref{implicitmap} we find
\begin{eqnarray}\label{inequlatyforyn}
\|y_n-p\| &=& \|\widehat\alpha_n(x_n-p)+\widehat\beta_n(I^{k(n)}_{j(n)}x_n-p)+\widehat\gamma_n(v_n-p)\|\nonumber\\
&\le& (1-\widehat\beta_n-\widehat\gamma_n)\|x_n-p\|+\widehat\beta_n\|I^{k(n)}_{j(n)}x_n-p\|+\widehat\gamma_n\|v_n-p\|\nonumber\\
&\le& (1-\widehat\beta_n-\widehat\gamma_n)\|x_n-p\|+\widehat\beta_n\mu_{k(n)}\|x_n-p\|+\widehat\gamma_n\|v_n-p\|\nonumber\\
&\le& (1-\widehat\beta_n)\|x_n-p\|+\widehat\beta_n\mu_{k(n)}\|x_n-p\|+\widehat\gamma_n\|v_n-p\|\nonumber\\
&\le& (1-\widehat\beta_n)\mu_{k(n)}\|x_n-p\|+\widehat\beta_n\mu_{k(n)}\|x_n-p\|+\widehat\gamma_n\|v_n-p\|\nonumber\\
&\le& \mu_{k(n)}\|x_n-p\|+\widehat\gamma_n\|v_n-p\| \nonumber\\
&\le& \lambda_{k(n)}\mu_{k(n)}\|x_n-p\|+\widehat\gamma_n\|v_n-p\|.
\end{eqnarray}
Then from \eqref{inequlatyforyn} we have
\begin{eqnarray*}
\|x_n-p\| &\le& (1-\beta_n)\|x_{n-1}-p\|+\beta_n\lambda^2_{k(n)}\mu^2_{k(n)}\|x_n-p\|\\
&&\qquad \qquad
+\gamma_n\|u_n-p\|+\beta_n\lambda_{k(n)}\mu_{k(n)}\widehat\gamma_n\|v_n-p\|,
\end{eqnarray*}
so one gets
\begin{eqnarray}\label{inequlatyforxn}
(1-\beta_n\lambda^2_{k(n)}\mu^2_{k(n)})\|x_n-p\| &\le&  (1-\beta_n)\|x_{n-1}-p\| \nonumber\\
&&+\gamma_n\|u_n-p\|+\beta_n\lambda_{k(n)}\mu_{k(n)}\widehat\gamma_n\|v_n-p\|.
\end{eqnarray}
By condition (iii) we obtain $\beta_n\lambda^2_{k(n)}\mu^2_{k(n)}\le
B^{*}\Lambda^2 M^2<1,$ and therefore
$$1-\beta_n\lambda^2_{k(n)}\mu^2_{k(n)}\ge 1-B^{*}\Lambda^2 M^2>0.$$ Hence
from \eqref{inequlatyforxn} we derive
\begin{eqnarray*}
\|x_n-p\| &\le&
\frac{1-\beta_n}{1-\beta_n\lambda^2_{k(n)}\mu^2_{k(n)}}\|x_{n-1}-p\|\\
&&\qquad\qquad+\frac{\gamma_n\|u_n-p\|+\beta_n\lambda_{k(n)}\mu_{k(n)}\widehat\gamma_n\|v_n-p\|}{1-\beta_n\lambda^2_{k(n)}\mu^2_{k(n)}}\\
&=&
\left(1+\frac{(\lambda^2_{k(n)}\mu^2_{k(n)}-1)\beta_n}{1-\beta_n\lambda^2_{k(n)}\mu^2_{k(n)}}\right)\|x_{n-1}-p\|\\
&&\qquad\qquad +\frac{\gamma_n\|u_n-p\|+\beta_n\lambda_{k(n)}\mu_{k(n)}\widehat\gamma_n\|v_n-p\|}{1-\beta_n\lambda^2_{k(n)}\mu^2_{k(n)}}\\
&\le&
\left(1+\frac{(\lambda^2_{k(n)}\mu^2_{k(n)}-1)\beta_n}{1-B^{*}\Lambda^2M^2}\right)\|x_{n-1}-p\|\\
&&\qquad\qquad
+\frac{\gamma_n\|u_n-p\|+\beta_n\lambda_{k(n)}\mu_{k(n)}\widehat\gamma_n\|v_n-p\|}{1-B^{*}\Lambda^2M^2}.
\end{eqnarray*}
Let
$$b_n=\dfrac{(\lambda^2_{k(n)}\mu^2_{k(n)}-1)\beta_n}{1-B^{*}\Lambda^2M^2}, \
\ \ \ \
c_n=\dfrac{\gamma_n\|u_n-p\|+\beta_n\lambda_{k(n)}\mu_{k(n)}\widehat\gamma_n\|v_n-p\|}{1-B^{*}\Lambda^2M^2}.$$
Then the last inequality has the following form
\begin{eqnarray}\label{VIinequalityforxn}
\|x_n-p\| &\le& \left(1+b_n\right)\|x_{n-1}-p\|+c_n.
\end{eqnarray}
From the condition (ii) we find
\begin{eqnarray*}
\sum\limits_{n=1}^\infty b_n &=&
\frac{1}{1-B^{*}\Lambda^2M^2}\sum\limits_{n=1}^{\infty}(\lambda^2_{k(n)}\mu^2_{k(n)}-1)\beta_n\\
&=&\frac{1}{1-B^{*}\Lambda^2M^2}\sum\limits_{n=1}^{\infty}(\lambda_{k(n)}\mu_{k(n)}-1)(\lambda_{k(n)}\mu_{k(n)}+1)\beta_n\\
&\le& \frac{\Lambda
M+1}{1-B^{*}\Lambda^2M^2}\sum\limits_{n=1}^{\infty}(\lambda_{k(n)}\mu_{k(n)}-1)\beta_n<\infty,
\end{eqnarray*}
and boundedness of the sequences $\{\|u_n-p\|\},$ $\{\|v_n-p\|\}$
with (iv) implies
\begin{eqnarray*}
\sum\limits_{n=1}^\infty c_n &=&
\sum\limits_{n=1}^\infty\frac{\gamma_n\|u_n-p\|+\beta_n\lambda_{k(n)}\mu_{k(n)}\widehat\gamma_n\|v_n-p\|}{1-B^{*}\Lambda^2M^2}\\
&\le&
\frac{1}{1-B^{*}\Lambda^2M^2}\sum\limits_{n=1}^\infty\gamma_n\|u_n-p\|+\frac{B^{*}\Lambda
M}{1-B^{*}\Lambda^2M^2}\sum\limits_{n=1}^\infty\widehat\gamma_n\|v_n-p\|<\infty.
\end{eqnarray*}
Now taking $a_n=\|x_{n-1}-p\|$ in \eqref{VIinequalityforxn} we
obtain
\begin{eqnarray*} a_{n+1}\le(1+b_n)a_n+c_n,
\end{eqnarray*}
and according to Lemma \ref{convergean} the limit
$\lim\limits_{n\to\infty}a_n$ exists. This means the limit
\begin{eqnarray}\label{limofxnminusp}
\lim\limits_{n\to\infty}\|x_n-p\|=d
\end{eqnarray}
exists, where $d\ge0$ is a constant. This completes the proof.
\end{pf}

Now we are ready to prove a general criteria of strong convergence
of \eqref{implicitmap}.

\begin{theorem}\label{criteria}
Let $X$ be a real Banach space and $K$ be a nonempty closed convex
subset of $X.$ Let $\{T_j\}_{j=1}^N:K\to K$ be a finite family of
uniformly $L_1-$Lipschitzian asymptotically quasi $I_j-$nonexpansive
mappings with a common sequence $\{\lambda_n\}\subset[1,\infty)$ and
$\{I_j\}_{j=1}^N:K\to K$ be a finite family of uniformly
$L_2-$Lipschitzian asymptotically quasi-nonexpansive mappings with a
common sequence $\{\mu_n\}\subset[1,\infty)$ such that
$F=\bigcap\limits_{j=1}^N \left(F(T_j)\cap
F(I_j)\right)\neq\emptyset.$  Suppose
$B^{*}=\sup\limits_{n}\beta_n,$
$\Lambda=\sup\limits_{n}\lambda_n\ge1,$ $M=\sup\limits_{n}\mu_n\ge1$
and $\{\alpha_n\},$ $\{\beta_n\},$ $\{\gamma_n\},$
$\{\widehat{\alpha}_n\},$ $\{\widehat{\beta}_n\},$
$\{\widehat{\gamma}_n\}$ are six sequences in $[0,1]$  which satisfy
the following conditions:
\begin{itemize}
  \item[(i)] $\alpha_n+\beta_n+\gamma_n=\widehat{\alpha}_n+\widehat{\beta}_n+\widehat{\gamma}_n=1, \ \ \ \forall n\ge1,$
  \item[(ii)] $\sum\limits_{n=1}^{\infty}(\lambda_n\mu_n-1)\beta_n<\infty,$
  \item[(iii)] $B^{*}<\dfrac{1}{\Lambda^2M^2},$
  \item[(iv)] $\sum\limits_{n=1}^\infty\gamma_n<\infty, \ \
  \sum\limits_{n=1}^\infty\widehat{\gamma}_n<\infty.$
\end{itemize}
Then the implicit iterative sequence $\{x_n\}$ with errors defined
by \eqref{implicitmap} converges strongly to a common fixed point in
$F$ if and only if
\begin{eqnarray}\label{conditionforcriteria}
\liminf\limits_{n\to\infty}d(x_n,F)=0.
\end{eqnarray}
\end{theorem}

\begin{pf}
The necessity of  condition \eqref{conditionforcriteria}  is
obvious. Let us proof the sufficiency part of the theorem.

Since $T_j,I_j:K\to K$ are uniformly $L_1,L_2-$Lipschitzian
mappings, respectively,  $T_j$ and $I_j$ are continuous mappings,
for each $j=\overline{1,N}$. Therefore, the sets $F(T_j)$ and
$F(I_j)$ are closed, for each $j=\overline{1,N}$. Hence
$F=\bigcap\limits_{j=1}^N \left(F(T_j)\cap F(I_j)\right)$ is a
nonempty closed set.

For any given $p\in F,$ we have (see \eqref{VIinequalityforxn})
\begin{eqnarray}\label{xnminusp}
\|x_n-p\| &\le& \left(1+b_n\right)\|x_{n-1}-p\|+c_n,
\end{eqnarray}
where $\sum\limits_{n=1}^\infty b_n<\infty$ and
$\sum\limits_{n=1}^\infty c_n<\infty.$  Hence, we find
\begin{eqnarray}\label{dxnF}
d(x_n,F) &\le& \left(1+b_n\right)d(x_{n-1},F)+c_n
\end{eqnarray}
So, the inequality \eqref{dxnF} with Lemma \ref{convergean} implies
the existence of the limit $\lim\limits_{n\to\infty}d(x_n,F)$. By
condition \eqref{conditionforcriteria}, one gets
\begin{eqnarray*}
\lim\limits_{n\to\infty}d(x_n,F)&=&
\liminf\limits_{n\to\infty}d(x_n,F)=0.
\end{eqnarray*}

Let us prove that the sequence $\{x_n\}$ converges to a common fixed
point of  $\{T_j\}_{j=1}^N$ and  $\{I_j\}_{j=1}^N.$

In fact, due to $1+t\le \exp(t)$ for all $t>0,$ and from
\eqref{xnminusp}, one finds
\begin{eqnarray}\label{xnminuspandexp}
\|x_n-p\| &\le& \exp(b_n)\|x_{n-1}-p\|+c_n.
\end{eqnarray}
Hence, for any positive integer $m,n,$ from \eqref{xnminuspandexp}
we find
\begin{eqnarray}\label{xnplusmp}
\|x_{n+m}-p\| &\le& \exp(b_{n+m})\|x_{n+m-1}-p\|+c_{n+m}\nonumber\\
&\le& \exp(b_{n+m}+b_{n+m-1})\|x_{n+m-2}-p\|\nonumber\\
&&+c_{n+m}+c_{n+m-1}\exp(b_{n+m})\nonumber\\
&\le& \exp(b_{n+m}+b_{n+m-1}+b_{n+m-2})\|x_{n+m-3}-p\|\nonumber\\
&&+c_{n+m}+c_{n+m-1}\exp(b_{n+m})+c_{n+m-2}\exp(b_{n+m}+b_{n+m-1})\nonumber\\
&&\vdots\nonumber \\
&\le& \exp\left(\sum\limits_{i=n+1}^{n+m} b_i\right)\|x_n-p\|+c_{n+m}+\sum\limits_{j=n+1}^{n+m-1}c_{j}\exp\left(\sum\limits_{i=j+1}^{n+m} b_i\right)\nonumber\\
&\le& \exp\left(\sum\limits_{i=n+1}^{n+m} b_i\right)\|x_n-p\|\nonumber\\
&&+c_{n+m}\exp\left(\sum\limits_{i=n+1}^{n+m}
b_i\right)+\sum\limits_{j=n+1}^{n+m-1}c_{j}\exp\left(\sum\limits_{i=n+1}^{n+m}
b_i\right)\nonumber\\
&\le& \exp\left(\sum\limits_{i=n+1}^{n+m}
b_i\right)\left(\|x_n-p\|+\sum\limits_{j=n+1}^{n+m}c_j\right)\nonumber\\
&\le& \exp\left(\sum\limits_{i=1}^{\infty}
b_i\right)\left(\|x_n-p\|+\sum\limits_{j=n+1}^{\infty}c_j\right)\nonumber\\
&\le& W\left(\|x_n-p\|+\sum\limits_{j=n+1}^{\infty}c_j\right),
\end{eqnarray}
for all $p\in F$, where $W=\exp\left(\sum\limits_{i=1}^{\infty}
b_i\right)<\infty.$

Since $\lim\limits_{n\to\infty}d(x_n,F)=0$ and
$\sum\limits_{j=1}^{\infty} c_j<\infty,$  for any given $\varepsilon
>0,$ there exists a positive integer number $n_0$ such that
\begin{eqnarray*}d(x_{n_0},F)<\frac{\varepsilon}{2W}, \qquad
\sum\limits_{j=n_0+1}^{\infty} c_j<\frac{\varepsilon}{2W}.
\end{eqnarray*}
Therefore there exists $p_1\in F$ such that
\begin{eqnarray*}\|x_{n_0}-p_1\|<\frac{\varepsilon}{2W},\qquad
\sum\limits_{j=n_0+1}^{\infty} c_j<\frac{\varepsilon}{2W}.
\end{eqnarray*}
Consequently, for all $n\ge n_0$  from \eqref{xnplusmp} we have
\begin{eqnarray*}
\|x_{n}-p_1\| &\le& W\left(\|x_{n_0}-p_1\|+\sum\limits_{j=n_0+1}^{\infty}c_j\right)\\
&<&W\cdot\frac{\varepsilon}{2W}+W\cdot\frac{\varepsilon}{2W}\\
&=& \varepsilon,
\end{eqnarray*}
this means that the sequence $\{x_n\}$ converges strongly to a
common fixed point $p_1$ of $\{T_j\}_{j=1}^N$ and $\{I_j\}_{j=1}^N.$
This completes the proof.
\end{pf}

To prove main results we need one more an auxiliary result.

\begin{proposition}\label{xnTxn&xnIxn}
Let $X$ be a real uniformly convex Banach space and $K$ be a
nonempty closed convex subset of $X.$ Let $\{T_j\}_{j=1}^N:K\to K$
be a finite family of uniformly $L_1-$Lipschitzian asymptotically
quasi $I_j-$nonexpansive mappings with a common sequence
$\{\lambda_n\}\subset[1,\infty)$ and $\{I_j\}_{j=1}^N:K\to K$ be a
finite family of uniformly $L_2-$Lipschitzian asymptotically
quasi-nonexpansive mappings with a common sequence
$\{\mu_n\}\subset[1,\infty)$ such that $F=\bigcap\limits_{j=1}^N
\left(F(T_j)\cap F(I_j)\right)\neq\emptyset.$ Suppose
$B_{*}=\inf\limits_{n}\beta_n,$ $B^{*}=\sup\limits_{n}\beta_n,$
$\Lambda=\sup\limits_{n}\lambda_n\ge1,$ $M=\sup\limits_{n}\mu_n\ge1$
and $\{\alpha_n\},$ $\{\beta_n\},$ $\{\gamma_n\},$
$\{\widehat{\alpha}_n\},$ $\{\widehat{\beta}_n\},$
$\{\widehat{\gamma}_n\}$ are six sequences in $[0,1]$  which satisfy
the following conditions:
\begin{itemize}
  \item[(i)] $\alpha_n+\beta_n+\gamma_n=\widehat{\alpha}_n+\widehat{\beta}_n+\widehat{\gamma}_n=1, \ \ \ \forall n\ge1,$
  \item[(ii)] $\sum\limits_{n=1}^{\infty}(\lambda_n\mu_n-1)\beta_n<\infty,$
  \item[(iii)] $0<B_{*}\le B^{*}<\dfrac{1}{\Lambda^2M^2}<1,$
  \item[(iv)] $0<\widehat B_{*}=\inf\limits_{n}\widehat\beta_n\le\sup\limits_{n}\widehat\beta_n=\widehat B^{*}<1,$
  \item[(v)] $\sum\limits_{n=1}^\infty\gamma_n<\infty, \ \
  \sum\limits_{n=1}^\infty\widehat{\gamma}_n<\infty.$
\end{itemize}
Then the implicit iterative sequence $\{x_n\}$ with errors defined
by \eqref{implicitmap}  satisfies the following
$$\lim\limits_{n\to\infty}\|x_n-T_jx_n\|=0, \qquad \lim\limits_{n\to\infty}\|x_n-I_jx_n\|=0, \qquad \forall j=\overline{1,N}.$$
\end{proposition}

\begin{pf} First, we shall prove that
$$\lim\limits_{n\to\infty}\|x_n-T^{k(n)}_{j(n)}x_n\|=0, \qquad \lim\limits_{n\to\infty}\|x_n-I^{k(n)}_{j(n)}x_n\|=0.$$

According to Lemma \ref{limexistsxnminusp} for any $p\in F$ we have
\begin{eqnarray}\label{xnminuspequald}
\lim\limits_{n\to\infty}\|x_n-p\|=d.
\end{eqnarray}
So, the sequence $\{x_n\}$ is bounded in $K$.

It follows from \eqref{implicitmap} that
\begin{eqnarray}\label{xnpnd}
\|x_n-p\|&=&
\|(1-\beta_n)(x_{n-1}-p+\gamma_n(u_n-x_{n-1}))\nonumber\\
&&\qquad
\qquad+\beta_n(T^{k(n)}_{j(n)}y_n-p+\gamma_n(u_n-x_{n-1}))\|.
\end{eqnarray}
Due to condition (v) and boundedness of the sequences $\{u_n\}$ and
$\{x_n\}$ we have
\begin{eqnarray}\label{xnpgammanun}
\limsup\limits_{n\to\infty}\|x_{n-1}&-&p+\gamma_n(u_n-x_{n-1})\|\le \nonumber \\
&\le&
\limsup\limits_{n\to\infty}\|x_{n-1}-p\|+\limsup\limits_{n\to\infty}\gamma_n\|u_n-x_{n-1}\|=d.\qquad\quad
\end{eqnarray}
By means of  asymptotically quasi $I_j-$nonexpansivity of $T_j$ and
asymptotically quasi-nonexpansivity of $I_j$ from
\eqref{inequlatyforyn} and boundedness of $\{u_n\},$ $\{v_n\},$
$\{x_n\}$ with condition (v) we obtain
\begin{eqnarray}\label{Tnynp}
\limsup\limits_{n\to\infty}&\|&T^{k(n)}_{j(n)}y_n-p+\gamma_n(u_n-x_{n-1})\ \ \ \|\nonumber\\
&\le&
\limsup\limits_{n\to\infty}\lambda_{k(n)}\mu_{k(n)}\|y_n-p\|+\limsup\limits_{n\to\infty}\gamma_n\|u_n-x_{n-1}\|\nonumber\\
&\le& \limsup\limits_{n\to\infty}\|y_n-p\|\nonumber\\
&\le&
\limsup\limits_{n\to\infty}\lambda_{k(n)}\mu_{k(n)}\|x_n-p\|+\limsup\limits_{n\to\infty}\widehat\gamma_n\|v_n-p\|=d.
\end{eqnarray}
Now using \eqref{xnpgammanun}, \eqref{Tnynp} and applying to Lemma
\ref{convexxnyn} to \eqref{xnpnd} one finds
\begin{eqnarray}\label{xnTnyn}
\lim\limits_{n\to\infty}\|x_{n-1}-T^{k(n)}_{j(n)}y_n\|=0.
\end{eqnarray}
From \eqref{implicitmap}, \eqref{xnTnyn} and condition (v) we infer
that
\begin{eqnarray}\label{xnxnminus1}
\lim\limits_{n\to\infty}\|x_n-x_{n-1}\| &=&
\lim\limits_{n\to\infty}\|\beta_n(T^{k(n)}_{j(n)}y_n-x_{n-1})+\gamma_n(u_n-x_{n-1})\|=0.\qquad
\end{eqnarray}
From \eqref{xnxnminus1} one can get
\begin{eqnarray}\label{xnxnminusj}
\lim\limits_{n\to\infty}\|x_n-x_{n+j}\| &=& 0, \qquad
j=\overline{1,N}.
\end{eqnarray}
On the other hand, we have
\begin{eqnarray*}
\|x_{n-1}-p\| &\le&
\|x_{n-1}-T^{k(n)}_{j(n)}y_n\|+\|T^{k(n)}_{j(n)}y_n-p\|\\
&\le&
\|x_{n-1}-T^{k(n)}_{j(n)}y_n\|+\lambda_{k(n)}\mu_{k(n)}\|y_n-p\|,
\end{eqnarray*}
which means
\begin{eqnarray*}
\|x_{n-1}-p\|-\|x_{n-1}-T^{k(n)}_{j(n)}y_n\|\le\lambda_{k(n)}\mu_{k(n)}\|y_n-p\|.
\end{eqnarray*}
The last inequality with \eqref{inequlatyforyn} implies that
\begin{eqnarray*}
\|x_{n-1}-p\|-\|x_{n-1}-T^{k(n)}_{j(n)}y_n\|&\le&
\lambda_{k(n)}\mu_{k(n)}\|y_n-p\|\le\\
&\le&
\lambda^2_{k(n)}\mu^2_{k(n)}\|x_n-p\|+\lambda_{k(n)}\mu_{k(n)}\widehat\gamma_n\|v_n-p\|.
\end{eqnarray*}
Then condition (v) and \eqref{xnTnyn}, \eqref{xnminuspequald} with
the Squeeze theorem yield
\begin{eqnarray}\label{ynpd}
\lim\limits_{n\to\infty}\|y_n-p\|=d
\end{eqnarray}
Again from \eqref{implicitmap} we can see that
\begin{eqnarray}\label{ynptod}
\|y_n-p\|=\|(1-\widehat\beta_n)(x_n-p&+&\widehat\gamma_n(v_n-x_n))\nonumber\\
&+&\widehat\beta_n(I^{k(n)}_{j(n)}x_n-p+\widehat\gamma_n(v_n-x_n))\|.
\end{eqnarray}
From  \eqref{xnminuspequald} with condition (v) one finds
\begin{eqnarray*}
\limsup\limits_{n\to\infty}\|x_n&-&p+\widehat\gamma_n(v_n-x_n)\|\le \nonumber \\
&\le&
\limsup\limits_{n\to\infty}\|x_n-p\|+\limsup\limits_{n\to\infty}\widehat\gamma_n\|v_n-x_n\|=d.\qquad\quad
\end{eqnarray*}
and
\begin{eqnarray*}
\limsup\limits_{n\to\infty}\|I^{k(n)}_{j(n)}x_n&-&p+\widehat\gamma_n(v_n-x_n)\|\nonumber\\
&\le&
\limsup\limits_{n\to\infty}\mu_{k(n)}\|x_n-p\|+\limsup\limits_{n\to\infty}\widehat\gamma_n\|v_n-x_{n}\|\nonumber\\
&\le& \limsup\limits_{n\to\infty}\|x_n-p\|=d.
\end{eqnarray*}
Now applying Lemma \ref{convexxnyn} to \eqref{ynptod} we obtain
\begin{eqnarray}\label{xnInxn}
\lim\limits_{n\to\infty}\|x_n-I^{k(n)}_{j(n)}x_n\|=0.
\end{eqnarray}
Consider
\begin{eqnarray*}
\|x_n-T^{k(n)}_{j(n)}x_n\| &\le&
\|x_n-x_{n-1}\|+\|x_{n-1}-T^{k(n)}_{j(n)}y_n\|+\|T^{k(n)}_{j(n)}y_n-T^{k(n)}_{j(n)}x_n\|\\
&\le& \|x_n-x_{n-1}\|+\|x_{n-1}-T^{k(n)}_{j(n)}y_n\|+L_1\|y_n-x_n\|\\
&=& \|x_n-x_{n-1}\|+\|x_{n-1}-T^{k(n)}_{j(n)}y_n\|\\
&&+L_1\|\widehat\beta_n(I^{k(n)}_{j(n)}x_n-x_n)+\widehat\gamma_n(v_n-x_n)\|\\
&\le& \|x_n-x_{n-1}\|+\|x_{n-1}-T^{k(n)}_{j(n)}y_n\|\\
&&+L_1\widehat\beta_n\|I^{k(n)}_{j(n)}x_n-x_n\|+L_1\widehat\gamma_n\|v_n-x_n\|.
\end{eqnarray*}
Then from \eqref{xnTnyn}, \eqref{xnxnminus1}, \eqref{xnInxn} and
condition (v) we get
\begin{eqnarray}\label{xnTnxn}
\lim\limits_{n\to\infty}\|x_n-T^{k(n)}_{j(n)}x_n\|=0.
\end{eqnarray}
Now we prove that
$$\lim\limits_{n\to\infty}\|x_n-T_jx_n\|=0, \qquad \lim\limits_{n\to\infty}\|x_n-I_jx_n\|=0, \qquad \forall j=\overline{1,N}.$$

For each $n>N$ we have $n\equiv n-N({\rm{mod}} N)$ and
$n=(k(n)-1)N+j(n),$ hence
$n-N=((k(n)-1)-1)N+j(n)=(k(n-N)-1)N+j(n-N),$ i.e.
$$k(n-N)=k(n)-1, \qquad j(n-N)=j(n).$$
So letting $T_n :=T_{j(n)({\rm{mod}} N)}$ we obtain
\begin{eqnarray*}
\|x_n-T_nx_n\| &\le& \|x_n-T^{k(n)}_{j(n)}x_n\|+\|T^{k(n)}_{j(n)}x_n-T_{j(n)}x_n\|\\
&\le& \|x_n-T^{k(n)}_{j(n)}x_n\| +L_1\|T^{{k(n)}-1}_{j(n)}x_n-x_n\|\\
&\le& \|x_n-T^{k(n)}_{j(n)}x_n\|+ L_1\|T^{{k(n)}-1}_{j(n)}x_n-T^{{k(n)}-1}_{j(n-N)}x_{n-N}\|\\
&& +L_1\|T^{{k(n)}-1}_{j(n-N)}x_{n-N}-x_{n-N}\|+L_1\|x_{n-N}-x_n\|\\
&\le& \|x_n-T^{k(n)}_{j(n)}x_n\| +L^2_1\|x_n-x_{n-N}\| \\
&& +L_1\|T^{{k(n-N)}}_{j(n-N)}x_{n-N}-x_{n-N}\|+L_1\|x_{n-N}-x_n\|\\
&\le& \|x_n-T^{k(n)}_{j(n)}x_n\|
+L_1(L_1+1)\|x_n-x_{n-N}\|\\
&&+L_1\|T^{{k(n-N)}}_{j(n-N)}x_{n-N}-x_{n-N}\|
\end{eqnarray*}
which with \eqref{xnxnminusj}, \eqref{xnTnxn} implies
\begin{eqnarray}\label{xnTxn}
\lim\limits_{n\to\infty}\|x_n-T_nx_n\|=0.
\end{eqnarray}
Analogously, one has
\begin{eqnarray*}
\|x_n-I_nx_n\| &\le& \|x_n-I^{k(n)}_{j(n)}x_n\|
+L_2(L_2+1)\|x_n-x_{n-N}\|\\
&&+L_2\|I^{{k(n-N)}}_{k(n-N)}x_{n-N}-x_{n-N}\|.
\end{eqnarray*}
which with \eqref{xnxnminusj}, \eqref{xnInxn} implies
\begin{eqnarray}\label{xnIxn}
\lim\limits_{n\to\infty}\|x_n-I_nx_n\|=0.
\end{eqnarray}
For any $j=\overline{1,N},$ from \eqref{xnxnminusj} and
\eqref{xnTxn} we have
\begin{eqnarray*}
\|x_n-T_{n+j}x_n\|&\le&
\|x_n-x_{n+j}\|+\|x_{n+j}-T_{n+j}x_{n+j}\|+\|T_{n+j}x_{n+j}-T_{n+j}x_n\|\\
&\le& (1+L_1)\|x_n-x_{n+j}\|+\|x_{n+j}-T_{n+j}x_{n+j}\|\to 0 \ \
(n\to\infty),
\end{eqnarray*}
which implies that the sequence
\begin{eqnarray}\label{xnTnjxn}
\bigcup\limits_{j=1}^N\bigl\{\|x_n-T_{n+j}x_n\|\bigr\}_{n=1}^\infty\to
0 \ \ (n\to\infty).
\end{eqnarray}
Analogously we have
\begin{eqnarray*}
\|x_n-I_{n+j}x_n\| &\le&
(1+L_2)\|x_n-x_{n+j}\|+\|x_{n+j}-I_{n+j}x_{n+j}\|\to 0 \ \
(n\to\infty),
\end{eqnarray*}
and
\begin{eqnarray}\label{xnInjxn}
\bigcup\limits_{j=1}^N\bigl\{\|x_n-I_{n+j}x_n\|\bigr\}_{n=1}^\infty\to
0 \ \ (n\to\infty).
\end{eqnarray}
According to
\begin{eqnarray*}
\bigl\{\|x_n-T_jx_n\|\bigr\}_{n=1}^\infty &=&
\bigl\{\|x_n-T_{n+(j-n)}x_n\|\bigr\}_{n=1}^\infty\\
&=& \bigl\{\|x_n-T_{n+j_n}x_n\|\bigr\}_{n=1}^\infty\subset
\bigcup\limits_{l=1}^N\bigl\{\|x_n-T_{n+l}x_n\|\bigr\}_{n=1}^\infty,
\end{eqnarray*}
and
\begin{eqnarray*}
\bigl\{\|x_n-I_jx_n\|\bigr\}_{n=1}^\infty &=&
\bigl\{\|x_n-I_{n+(j-n)}x_n\|\bigr\}_{n=1}^\infty\\
&=& \bigl\{\|x_n-I_{n+j_n}x_n\|\bigr\}_{n=1}^\infty\subset
\bigcup\limits_{l=1}^N\bigl\{\|x_n-I_{n+l}x_n\|\bigr\}_{n=1}^\infty,
\end{eqnarray*}
where $j-n\equiv j_n({\rm{mod}} N),$ $j_n\in \{1,2,\cdots,N\},$
 from \eqref{xnTnjxn}, \eqref{xnInjxn} we find
\begin{eqnarray}\label{xnTjxnxnIjxn}
\lim\limits_{n\to\infty}\|x_n-T_jx_n\|=0, \qquad
\lim\limits_{n\to\infty}\|x_n-I_jx_n\|=0, \qquad \forall
j=\overline{1,N}.
\end{eqnarray}
\end{pf}

Now we are ready to formulate one of main result concerning weak
convergence of the sequence $\{x_n\}$.

\begin{theorem}\label{weakconvergence}
Let $X$ be a real uniformly convex Banach space satisfying Opial
condition and $K$ be a nonempty closed convex subset of $X.$ Let
$E:X\to X$ be an identity mapping, $\{T_j\}_{j=1}^N:K\to K$ be a
finite family of uniformly $L_1-$Lipschitzian asymptotically quasi
$I_j-$nonexpansive mappings with a common sequence
$\{\lambda_n\}\subset[1,\infty)$ and $\{I_j\}_{j=1}^N:K\to K$ be a
finite family of uniformly $L_2-$Lipschitzian asymptotically
quasi-nonexpansive mappings with a common sequence
$\{\mu_n\}\subset[1,\infty)$ such that $F=\bigcap\limits_{j=1}^N
\left(F(T_j)\cap F(I_j)\right)\neq\emptyset.$ Suppose
$B_{*}=\inf\limits_{n}\beta_n,$ $B^{*}=\sup\limits_{n}\beta_n,$
$\Lambda=\sup\limits_{n}\lambda_n\ge1,$ $M=\sup\limits_{n}\mu_n\ge1$
and $\{\alpha_n\},$ $\{\beta_n\},$ $\{\gamma_n\},$
$\{\widehat{\alpha}_n\},$ $\{\widehat{\beta}_n\},$
$\{\widehat{\gamma}_n\}$ are six sequences in $[0,1]$  which satisfy
the following conditions:
\begin{itemize}
  \item[(i)] $\alpha_n+\beta_n+\gamma_n=\widehat{\alpha}_n+\widehat{\beta}_n+\widehat{\gamma}_n=1, \ \ \ \forall n\ge1,$
  \item[(ii)] $\sum\limits_{n=1}^{\infty}(\lambda_n\mu_n-1)\beta_n<\infty,$
  \item[(iii)] $0<B_{*}\le B^{*}<\dfrac{1}{\Lambda^2M^2}<1,$
  \item[(iv)] $0<\widehat B_{*}=\inf\limits_{n}\widehat\beta_n\le\sup\limits_{n}\widehat\beta_n=\widehat B^{*}<1,$
  \item[(v)] $\sum\limits_{n=1}^\infty\gamma_n<\infty, \ \
  \sum\limits_{n=1}^\infty\widehat{\gamma}_n<\infty.$
\end{itemize}
If the mappings $E-T_j$ and $E-I_j$ are semi-closed at zero for
every $j=\overline{1,N},$ then the implicit iterative sequence
$\{x_n\}$ with errors defined by \eqref{implicitmap} converges
weakly to a common fixed point of finite families of asymptotically
quasi $I_j-$nonexpansive mappings $\{T_j\}_{j=1}^N$ and
asymptotically quasi-nonexpansive mappings $\{I_j\}_{j=1}^N.$
\end{theorem}
\begin{pf}
Let $p\in F$, the according to Lemma \ref{limexistsxnminusp} the
sequence $\{\|x_n-p\|\}$ converges. This provides that $\{x_n\}$ is
bounded. Since $X$ is uniformly convex, then every bounded subset of
$X$ is weakly compact. From boundedness of $\{x_n\}$ in $K,$ we
find a subsequence $\{x_{n_k}\}\subset\{x_n\}$ such that
$\{x_{n_k}\}$ converges weakly to $q\in K.$ Hence from
\eqref{xnTjxnxnIjxn}, it follows that
$$\lim\limits_{n_k\to\infty}\|x_{n_k}-T_jx_{n_k}\|=0, \  \lim\limits_{n_k\to\infty}\|x_{n_k}-I_jx_{n_k}\|=0, \qquad \forall j=\overline{1,N}.$$
Since the mappings $E-T_j$ and $E-I_j$ are semi-closed at zero,
therefore we have $T_jq=q$ and $I_jq=q,$ for all $j=\overline{1,N},$
which means $q\in F.$

Finally, we prove that $\{x_n\}$ converges weakly to $q.$ In fact,
suppose the contrary, then there exists some subsequence
$\{x_{n_j}\}\subset\{x_n\}$ such that $\{x_{n_j}\}$ converges weakly
to $q_1\in K$ and $q_1\neq q$. Then by the same method as given
above, we can also prove that $q_1\in F.$

Taking $p=q$ and $p=q_1$ and using the same argument given in the
proof of \eqref{limofxnminusp}, we can prove that the limits
$\lim\limits_{n\to\infty}\|x_n-q\|$ and
$\lim\limits_{n\to\infty}\|x_n-q_1\|$ exist, and we have
$$\lim\limits_{n\to\infty}\|x_n-q\|=d, \ \ \
\lim\limits_{n\to\infty}\|x_n-q_1\|=d_1,$$ where $d, d_1$ are two
nonnegative numbers. By virtue of the Opial condition of $X$, one
finds
\begin{eqnarray*} d &=& \limsup\limits_{n_k\to\infty}\|x_{n_k}-q\| <
\limsup\limits_{n_k\to\infty}\|x_{n_k}-q_1\|=\\
&=& \limsup\limits_{n_j\to\infty}\|x_{n_j}-q_1\| <
\limsup\limits_{n_j\to\infty}\|x_{n_j}-q\| = d
\end{eqnarray*}
This is a contradiction. Hence $q_1=q.$ This implies that $\{x_n\}$
converges weakly to $q.$ This completes the proof of Theorem
\ref{weakconvergence}.
\end{pf}

Next, we prove strong convergence theorem

\begin{theorem}\label{strongconvergence}
Let $X$ be a real uniformly convex Banach space and $K$ be a
nonempty closed convex subset of $X.$ Let   $\{T_j\}_{j=1}^N:K\to K$
be a finite family of uniformly $L_1-$Lipschitzian asymptotically
quasi $I_j-$nonexpansive mappings with a common sequence
$\{\lambda_n\}\subset[1,\infty)$ and $\{I_j\}_{j=1}^N:K\to K$ be a
finite family of uniformly $L_2-$Lipschitzian asymptotically
quasi-nonexpansive mappings with a common sequence
$\{\mu_n\}\subset[1,\infty)$ such that $F=\bigcap\limits_{j=1}^N
\left(F(T_j)\cap F(I_j)\right)\neq\emptyset.$ Suppose
$B_{*}=\inf\limits_{n}\beta_n,$ $B^{*}=\sup\limits_{n}\beta_n,$
$\Lambda=\sup\limits_{n}\lambda_n\ge1,$ $M=\sup\limits_{n}\mu_n\ge1$
and $\{\alpha_n\},$ $\{\beta_n\},$ $\{\gamma_n\},$
$\{\widehat{\alpha}_n\},$ $\{\widehat{\beta}_n\},$
$\{\widehat{\gamma}_n\}$ are six sequences in $[0,1]$  which satisfy
the following conditions:
\begin{itemize}
  \item[(i)] $\alpha_n+\beta_n+\gamma_n=\widehat{\alpha}_n+\widehat{\beta}_n+\widehat{\gamma}_n=1, \ \ \ \forall n\ge1,$
  \item[(ii)] $\sum\limits_{n=1}^{\infty}(\lambda_n\mu_n-1)\beta_n<\infty,$
  \item[(iii)] $0<B_{*}\le B^{*}<\dfrac{1}{\Lambda^2M^2}<1,$
  \item[(iv)] $0<\widehat B_{*}=\inf\limits_{n}\widehat\beta_n\le\sup\limits_{n}\widehat\beta_n=\widehat B^{*}<1,$
  \item[(v)] $\sum\limits_{n=1}^\infty\gamma_n<\infty, \ \
  \sum\limits_{n=1}^\infty\widehat{\gamma}_n<\infty.$
\end{itemize}
If at least one mapping of the mappings $\{T_j,I_j\}_{j=1}^N$ is
semi-compact, then the implicit iterative sequence $\{x_n\}$ with
errors defined by \eqref{implicitmap} converges strongly to a common
fixed point of finite families of asymptotically quasi
$I_j-$nonexpansive mappings $\{T_j\}_{j=1}^N$ and asymptotically
quasi-nonexpansive mappings $\{I_j\}_{j=1}^N.$
\end{theorem}

\begin{pf}
Without any loss of generality, we may assume that $T_1$ is semi-compact.
This with \eqref{xnTjxnxnIjxn} means that there exists a subsequence
$\{x_{n_k}\}\subset\{x_n\}$ such that $x_{n_k}\to x^{*}$ strongly
and $x^{*}\in K.$ Since $T_j,I_j$ are continuous, then from
\eqref{xnTjxnxnIjxn}, for all $j=\overline{1,N}$ we find
$$\|x^{*}-T_jx^{*}\|=\lim\limits_{n_k\to\infty}\|x_{n_k}-T_jx_{n_k}\|=0, \ \ \ \ \|x^{*}-I_jx^{*}\|=\lim\limits_{n_k\to\infty}\|x_{n_k}-I_jx_{n_k}\|=0.$$
This shows that $x^{*}\in F.$ According to Lemma
\ref{limexistsxnminusp} the limit
$\lim\limits_{n\to\infty}\|x_n-x^{*}\|$ exists. Then
$$\lim\limits_{n\to\infty}\|x_n-x^{*}\|=\lim\limits_{n_k\to\infty}\|x_{n_k}-x^{*}\|=0,$$
which means $\{x_n\}$ converges to $x^{*}\in F.$ This completes the
proof.
\end{pf}

\begin{remark} If we take $\gamma_n = 0$, for all $n\in\bn$ and $I_j=E$, $j=1,2,\dots,N$ then the above theorem becomes Theorem 3.3 due to Sun \cite{[Sun]}. If $I_j=E$, $j=1,2,\dots,N$
and $\{T_j\}_{j=1}^N$ are asymptotically nonexpansive, then
we got main results of \cite{CTL,gou}.
\end{remark}

\section*{acknowledgments} The authors acknowledges the MOSTI
grant 01-01-08-SF0079.

\end{document}